\documentclass[12pt]{amsart}
\usepackage{amsmath,amsthm,latexsym,amscd,amsbsy,amssymb}
 \relax

\chardef\bslash=`\\ 

\makeatletter
\def\verbatim{\interlinepenalty\@M \@verbatim
 \leftskip\@totalleftmargin\advance\leftskip2pc
 \frenchspacing\@vobeyspaces \@xverbatim}
\makeatother
\makeatletter
 \def\dgt@k{\dg@DX=-3 \dg@DY=2 \dg@SIZE=3}
\makeatother

\makeatletter
 \def\dgt@kk{\dg@DX=3 \dg@DY=-1 \dg@SIZE=3}%
\makeatother

\theoremstyle{plain}
\newtheorem{thm}{Theorem}[section]
\newtheorem{cor}[thm]{Corollary}
\newtheorem{lem}[thm]{Lemma}

\newtheorem*{A}{Theorem 3.1}

\theoremstyle{definition}

\newtheorem{defin}[thm]{Definition}

\newcommand{\ed}{\operatorname{edd}}

\numberwithin{equation}{section}

\def\cov{{\rm cov}}

\def\dim{{\rm dim}}

\def\ed{\text{{\rm e-dim}}}

\def\pr{{\rm pr}}

\def\St{{\rm St}}
\def\wh{\widehat}
\def\wt{\widetilde}

\begin{document}


\title[Hurewicz theorem for extension dimension]
{Hurewicz theorem for extension dimension}
\author{N.~Brodsky}
\address{Department of Mathematics and Statistics,
University of Saskatche\-wan,
McLean Hall, 106 Wiggins Road, Saskatoon, SK, S7N 5E6, Canada}
\email{brodsky@math.usask.ca}
\author{A.~Chigogidze}
\address{Department of Mathematics and Statistics,
University of Saskatche\-wan,
McLean Hall, 106 Wiggins Road, Saskatoon, SK, S7N 5E6, Canada}
\email{chigogid@math.usask.ca}
\thanks{The second author was partially supported by NSERC research grant.}


\begin{abstract}{We prove a new selection theorem for multivalued mappings
of C-space. Using this theorem we prove extension dimensional version
of Hurewicz theorem for a closed mapping $f\colon X\to Y$ of $k$-space $X$
onto paracompact $C$-space $Y$: if for finite $CW$-complex $M$
we have $\ed Y\le [M]$ and for every point $y\in Y$ and
every compactum $Z$ with $\ed Z\le [M]$ we have
$\ed(f^{-1}(y)\times Z)\le [L]$
for some $CW$-complex $L$, then $\ed X\le [L]$.}
\end{abstract}

\maketitle
\markboth{N.~Brodsky, A.~Chigogidze}
{Hurewicz theorem for extension dimension}

\section{Introduction}\label{S:intro}

The classical Hurewicz theorem states that for a mapping of
finite-dimensional compacta $f\colon X\to Y$ we have
\[ \dim X\le \dim Y+\dim f ,\;\; \text{where}\;\; \dim f=\max\{\dim(f^{-1}(y)\mid y\in Y\} . \]
There are several approaches to extension dimensional generalization of
Hure\-wicz theorem~\cite{DRS},\cite{ChV},\cite{BCK},\cite{Le},\cite{LL},\cite{T}.

Using the idea from~\cite{ChV} we improve Theorem~7.6 from~\cite{BCK}:

\begin{A}
Let $f\colon X\to Y$ be a closed mapping of a $k$-space $X$ onto
paracompact $C$-space $Y$.
Suppose that $\ed Y\le [M]$ for a finite $CW$-complex $M$.
If for every point $y\in Y$ and for every compactum $Z$
with $\ed Z\le [M]$ we have $\ed(f^{-1}(y)\times Z)\le [L]$
for some $CW$-complex $L$, then $\ed X\le [L]$.
\end{A}

The notion of extension dimension was introduced by Dranishnikov~\cite{D}:
for a $CW$-complex $L$ a space $X$ is said to have {\it extension dimension}
$\le [L]$ (notation: $\ed X \le [L]$) if any mapping of its closed subspace
$A\subset X$ into $L$ admits an extension to the whole space $X$.

To prove Theorem~\ref{thmHur} we need an extension dimensional version of
Uspenskij's selection theorem~\cite{U}. In section~\ref{S:Selections} we
prove Theorem~\ref{thmsel} on selections of multivalued mappings of C-space.
Then Theorem~\ref{thmtop} helps us
to prove Theorem~\ref{thmusp} --- a needed version of Uspenskij's theorem.

Filtrations of multivalued maps are proved to be very useful
for construction of continuous selections~\cite{SB},~\cite{BCK}.
And we state our selection theorems in terms of filtrations.
Note that Valov~\cite{V} used filtrations to prove a selection
theorem for mappings of finite $C$-spaces.

Let us recall some definitions and introduce our notations.
A space $X$ is called a {\it $k$-space} if $U\subset X$
is open in $X$ whenever $U\cap C$ is relatively open in $C$
for every compact subset $C$ of $X$.
The {\it graph} of a multivalued mapping $F\colon X\to Y$
is the subset $\Gamma_F=\{(x,y)\in X\times Y\colon y\in F(x)\}$
of the product $X\times Y$.

We denote by $\cov X$ the collection of all coverings of the space $X$.
For a cover $\omega$ of a space $X$ and for a subset $A \subseteq X$ let
$\St(A,\omega)$ denote the star of the set $A$ with respect to $\omega$.
We say that a subset $A\subset X$ {\it refines} a cover $\omega\in \cov X$
if $A$ is contained in some element of $\omega$.
A covering $\omega'\in \cov X$ {\it strongly star refines} a covering
$\omega\in \cov X$ if for any element $W\in\omega'$ the set
$\St(W,\omega')$ refines $\omega$.

\begin{defin}
A topological space $X$ is called {\it $C$-space} if for each sequence
$\{\omega_i\}_{i\ge 1}$ of open covers of $X$, there is an open cover
$\Sigma$ of $X$ of the form $\cup_{i=1}^\infty \sigma_i$
such that for each $i\ge 1$, $\sigma_i$ is a
pairwise disjoint collection which refines $\omega_i$.
\end{defin}

If the space $X$ is paracompact, we can choose the cover $\Sigma$
to be locally finite and every collection $\sigma_i$ to be discrete.

\begin{defin}
A multivalued mapping $F\colon X\to Y$ is said to be
{\it strongly lower semicontinuous} (briefly, strongly l.s.c.)
if for any point $x\in X$ and any compact set $K\subset F(x)$ there exists
a neighborhood $V$ of $x$ such that $K\subset F(z)$ for every $z\in V$.
\end{defin}

\begin{defin}
Let $L$ be a $CW$-complex. A pair of spaces $V\subset U$ is said
to be {\it $[L]$-connected} (resp., {\it $[L]_c$-connected})
if for every paracompact space $X$ (resp., compact metric space $X$)
of extension dimension $\ed X\le [L]$
and for every closed subspace $A\subset X$ any mapping of $A$
into $V$ can be extended to a mapping of $X$ into $U$.
\end{defin}

An increasing\footnote{We consider only increasing filtrations
indexed by a segment of the integral series.}
sequence of subspaces $ Z_0\subset Z_1\subset\dots\subset Z $
is called a {\it filtration} of space $Z$.
A sequence of multivalued mappings $\{F_k\colon X\to Y\}$
is called a {\it filtration of multivalued mapping} $F\colon X\to Y$
if $\{F_k(x)\}$ is a filtration of $F(x)$ for any $x\in X$.

\begin{defin}
A filtration of multivalued mappings $\{G_i\colon X\to Y\}$
is said to be {\it fiberwise $[L]_c$-connected} if for any point $x\in X$
and any $i$ the pair $G_i(x)\subset G_{i+1}(x)$ is $[L]_c$-connected.
\end{defin}

\section{Selection theorems}\label{S:Selections}

The following notion of stably $[L]$-connected filtration
of multivalued mappings provides a key property of
the filtration for our construction of continuous selections.

\begin{defin}
A pair $F\subset H$ of multivalued mappings from $X$ to $Y$
is called {\it stably $[L]$-connected} if every point $x\in X$
has a neighborhood $O_x$ such that
the pair $F(O_x)\subset \cap_{z\in O_x} H(z)$ is $[L]$-connected.

We say that the pair $F\subset H$ is called {\it stably $[L]$-connected
with respect to a covering $\omega\in\cov X$}, if for any $W\in \omega$
the pair $F(W)\subset \cap_{x\in W} H(x)$ is $[L]$-connected.

A filtration $\{F_{i}\}$ of multivalued mappings is called
{\it stably $[L]$-connected} if every pair
$F_{i}\subset F_{i+1}$ is stably $[L]$-connected.
\end{defin}

Clearly, any stably $[L]$-connected pair of multivalued maps of a
space $X$ is stably $[L]$-connected with respect to some covering of $X$.

We denote by $Q$ the Hilbert cube.
We identify a space $Y$ with the subspace $Y\times \{0\}$ of the product
$Y\times Q$ and denote by $\pr_Y$ the projection of $Y\times Q$ onto $Y$.

\begin{defin}
For a subspace $Z\subset Y\times Q$ we say that $Y$
{\it projectively contains} $Z$.
We say that a multivalued mapping $F\colon X\to Y$
{\it projectively contains} a multivalued mapping $G\colon X\to Y\times Q$
if for any point $x\in X$ the set $\pr_Y\circ G(x)$ is contained in $F(x)$.
\end{defin}

\begin{lem} \label{lemmapairs}
Let $L$ be a finite $CW$-complex.
If a topological space $Y$ contains a compactum $K$ of extension dimension
$\ed K\le [L]$ such that the pair $K\subset Y$ is $[L]_c$-connected, then
$Y$ projectively contains a compactum $K'$ of extension dimension
$\ed K'\le [L]$ such that $K$ lies in $K'$ and
the pair $K\subset K'$ is $[L]$-connected.
\end{lem}

\begin{proof}
There exists $AE([L])$-compactum $K'$ of extension dimension
$\ed K'\le [L]$ containing the given compactum $K$~\cite{Ch97}.
Clearly, the pair $K\subset K'$ is $[L]$-connected.
Since $\ed K'\le [L]$, there exists a mapping $p\colon K'\to Y$
extending the inclusion of $K$ into $Y$.

It is easy to see that there exists a mapping $q\colon K'\to Q$
such that $q^{-1}(0)=K$ and $q$ is an embedding on $K'\setminus K$.
Now define an embedding $j\colon K'\to Y\times Q$ as $j=p\times q$.
Since $q^{-1}(0)=K$, the mapping $j$ coincide with $p$
on $K$ which is inclusion on $K$.
\end{proof}

\begin{defin}
We say that a filtration $F_0\subset F_1\subset\dots$ of multivalued mappings
from $X$ to $Y$ {\it projectively contains} a filtration
$G_0\subset G_1\subset\dots$ of multivalued mappings from $X$ to $Y\times Q$
if for any point $x\in X$ and any $n$ the set
$\pr_Y\circ G_n(x)$ is contained in $F_n(x)$.
\end{defin}

\begin{thm} \label{thmtop}
For a finite $CW$-complex $L$
any fiberwise $[L]_c$-connected filtration of strongly l.s.c. multivalued
mappings of paracompact space $X$ to a topological space $Y$ projectively
contains stably $[L]$-connected filtration of compact-valued mappings.
\end{thm}

\begin{proof}
For a given fiberwise $[L]_c$-connected filtration
$F_0\subset F_1\subset\dots$ of strongly l.s.c.
multivalued mappings we construct stably $[L]$-connected filtration
$G_0\subset G_1\subset\dots$ of compact-valued mappings
$G_n\colon X\to Y\times Q^n$ as follows:
successively for every $n\ge 0$ we construct a covering
$\omega_n=\{W^n_\lambda\}_{\lambda\in\Lambda_n}\in\cov X$
and a family of subcompacta $\{K^n_\lambda\}_{\lambda\in\Lambda_n}$
of $Y\times Q^n$, and define the mapping $G_n$ by the formula
\[  G_n(x)=\cup \{K^n_\lambda\mid x\in W^n_\lambda\}.  \]

First, we construct $G_0$, i.e. the covering $\omega_0$
and the family $\{K^0_\lambda\}_{\lambda\in\Lambda_0}$.
Since $F_0$ is strongly l.s.c., there exists a locally finite open covering
$\omega_{-1}=\{W^{-1}_\lambda\}_{\lambda\in\Lambda_{-1}}\in\cov X$
and a family $\{M^{-1}_\lambda\}_{\lambda\in\Lambda_{-1}}$
of points in $Y$ such that
$W^{-1}_\lambda\times M^{-1}_\lambda\subset\Gamma_{F_0}$
for any $\lambda\in\Lambda_{-1}$.
Denote by $H_0$ a multivalued mapping taking a point $x\in X$
to the set $H_0(x)=\cup \{M^{-1}_\lambda\mid x\in W^{-1}_\lambda\}$.
Note that $H_0(x)$ is contained in $F_0(x)$ and consists of
finitely many points.
By Lemma~\ref{lemmapairs} for any $x\in X$ there exists a compactum
$\wh H_0(x)\subset F_1(x)\times Q$ of extension dimension
$\ed \wh H_0(x)\le [L]$ such that the pair $H_0(x)\subset \wh H_0(x)$
is $[L]$-connected.
Since $F_1$ is strongly l.s.c., any point $x\in X$ has a neighborhood
$\mathcal O_0(x)$ such that the product $\mathcal O_0(x)\times \wh H_0(x)$
is contained in $\Gamma_{F_1}\times Q$.
Since $X$ is paracompact, we can choose neighborhoods $\mathcal O_0(x)$
in such a way that the covering $\mathcal O_0=\{\mathcal O_0(x)\}_{x\in X}$
strongly star refines $\omega_{-1}$.
Let $\omega_{0}=\{W^{0}_\lambda\}_{\lambda\in\Lambda_{0}}$
be a locally finite open cover of $X$ refining $\mathcal O_0$.
For every $\lambda\in\Lambda_0$ we fix a point $x_\lambda$ such that
$W^{0}_\lambda\subset \mathcal O_0(x_\lambda)$ and put
$M^0_\lambda=\wh H_0(x_\lambda)$.
For every $\lambda\in\Lambda_0$ we fix $\alpha(\lambda)\in\Lambda_{-1}$
such that $\St(W^{0}_\lambda,\mathcal O_0)\subset W^{-1}_{\alpha(\lambda)}$
and put $K^0_\lambda=M^{-1}_{\alpha(\lambda)}$.

Inductive step of our construction is similar to the first step.
Suppose that a covering
$\omega_{n-1}=\{W^{n-1}_\lambda\}_{\lambda\in\Lambda_{n-1}}\in\cov X$
and a family $\{M^{n-1}_\lambda\}_{\lambda\in\Lambda_{n-1}}$
of compacta in $Y\times Q^{n-1}$ are already constructed such that
$\ed M^{n-1}_\lambda\le [L]$ and the product
$W^{n-1}_\lambda\times M^{n-1}_\lambda$
is contained in $\Gamma_{F_n}\times Q^n$
for any $\lambda\in\Lambda_{n-1}$.
Denote by $H_n$ a multivalued mapping taking a point $x\in X$
to the compactum $H_n(x)=\cup \{M^{n-1}_\lambda\mid x\in W^{n-1}_\lambda\}$.
Note that $H_n(x)$ is contained in $F_n(x)\times Q^n$
and has extension dimension $\ed H_n(x)\le [L]$.
By Lemma~\ref{lemmapairs} for any $x\in X$ there exists a compactum
$\wh H_n(x)\subset F_{n+1}(x)\times Q^{n+1}$ of extension dimension
$\ed \wh H_n(x)\le [L]$ such that the pair
$H_n(x)\subset \wh H_n(x)$ is $[L]$-connected.
Since $F_{n+1}$ is strongly l.s.c., any point $x\in X$ has a neighborhood
$\mathcal O_n(x)$ such that the product $\mathcal O_n(x)\times \wh H_n(x)$
is contained in $\Gamma_{F_{n+1}}\times Q^{n+1}$.
Since $X$ is paracompact, we can choose neighborhoods $\mathcal O_n(x)$
in such a way that the covering $\mathcal O_n=\{\mathcal O_n(x)\}_{x\in X}$
strongly star refines $\omega_{n-1}$.
Let $\omega_{n}=\{W^{n}_\lambda\}_{\lambda\in\Lambda_{n}}$
be a locally finite open cover of $X$ refining $\mathcal O_n$.
For every $\lambda\in\Lambda_n$ we fix a point $x_\lambda$ such that
$W^{n}_\lambda\subset \mathcal O_n(x_\lambda)$ and put
$M^n_\lambda=\wh H_n(x_\lambda)$.
For every $\lambda\in\Lambda_n$ we fix $\alpha(\lambda)\in\Lambda_{n-1}$
such that $\St(W^{n}_\lambda,\mathcal O_n)\subset W^{n-1}_{\alpha(\lambda)}$
and put $K^n_\lambda=M^{n-1}_{\alpha(\lambda)}$.

To show that the pair $G_{n-1}\subset G_n$ is stably $[L]$-connected,
we prove that the pair
$G_{n-1}(W^{n}_\lambda)\subset \cap \{G_n(x)\mid x\in W^{n}_\lambda\}$
is $[L]$-connected for any $W^{n}_\lambda\in\omega_n$.
By the construction of $G_n$, the set $K^n_\lambda$
is contained in $\cap \{G_n(x)\mid x\in W^{n}_\lambda\}$.
We know that the pair
$H_{n-1}(x_{\alpha(\lambda)})\subset\wh H_{n-1}(x_{\alpha(\lambda)})=
M^{n-1}_{\alpha(\lambda)}=K^n_\lambda$ is $[L]$-connected.
Therefore it is enough to show the following inclusion:
\begin{multline*}
G_{n-1}(W^{n}_\lambda)=\\\bigcup\{K^{n-1}_\beta\mid
W^{n}_\lambda\cap W^{n-1}_\beta\ne\emptyset\} \subset
\cup\{M^{n-2}_\nu\mid x_{\alpha(\lambda)}\in W^{n-2}_\nu\}=
H_{n-1}(x_{\alpha(\lambda)})
\end{multline*}

\noindent which follows from the fact that
$W^{n}_\lambda\cap W^{n-1}_\beta\ne\emptyset$ implies
$x_{\alpha(\lambda)}\in W^{n-2}_{\alpha(\beta)}$
(note that $M^{n-2}_{\alpha(\beta)}=K^{n-1}_\beta$).
By the choice of $\alpha(\lambda)$ we have
$W^{n}_\lambda\subset \mathcal O_{n-1}(x_{\alpha(\lambda)})$.
Then $W^{n}_\lambda\cap W^{n-1}_\beta\ne\emptyset$ implies
$\mathcal O_{n-1}(x_{\alpha(\lambda)})\cap W^{n-1}_\beta\ne\emptyset$
and $x_{\alpha(\lambda)}\in \mathcal O_{n-1}(x_{\alpha(\lambda)})
\subset \St(W^{n-1}_\beta, \mathcal O_{n-1})\subset W^{n-2}_{\alpha(\beta)}$.
\end{proof}

\begin{defin}
For a space $Z$ a pair of spaces $V\subset U$ is said to be
{\it $Z$-connected} if for every closed subspace $A\subset Z$ any
mapping of $A$ into $V$ can be extended to a mapping of $Z$ into $U$.
\end{defin}

\begin{defin}
A pair $F\subset H$ of multivalued mappings from $X$ to $Y$
is called {\it stably $Z$-connected} if every point $x\in X$
has a neighborhood $O_x$ such that
the pair $F(O_x)\subset \cap_{z\in O_x} H(z)$ is $Z$-connected.

We say that the pair $F\subset H$ is called {\it stably $Z$-connected
with respect to a covering $\omega\in\cov X$}, if for any $W\in \omega$
the pair $F(W)\subset \cap_{x\in W} H(x)$ is $Z$-connected.

A filtration $\{F_{i}\}$ of multivalued mappings is called
{\it stably $Z$-connected} if every pair
$F_{i}\subset F_{i+1}$ is stably $Z$-connected.
\end{defin}

\begin{thm} \label{thmsel}
Let $F\colon X\to Y$ be a multivalued mapping of paracompact C-space $X$
to a topological space $Y$.
If $F$ admits infinite stably $X$-connected filtration of multivalued
mappings, then $F$ has a singlevalued continuous selection.
\end{thm}

\begin{proof}
Let $\{F_i\}_{i=-1}^\infty$ be the given filtration of $F$.
Let $\{\omega_i\}_{i=-1}^\infty$ be a sequence of coverings of $X$
such that $\omega_{i+1}$ refines $\omega_i$ and the pair
$F_{i}\subset F_{i+1}$ is stably $X$-connected
with respect to the covering $\omega_i$.
Since $X$ is paracompact C-space, there exists a locally finite
closed cover $\Sigma$ of $X$ of the form $\Sigma=\cup_{i=0}^\infty \sigma_i$
such that $\sigma_i$ is discrete collection refining $\omega_i$.
Define $\Sigma_n=\cup_{i=0}^n \sigma_i$.
We will construct a continuous selection $f$ of $F$
extending it successively over the sets $\Sigma_n$.

First, we construct $f_0\colon \Sigma_0\to Y$.
We define $f_0$ separately on every element $s$ of the discrete
collection $\sigma_0$: take a point $p\in F_{-1}(s)$ and put $f_0(s)=p$.
Since the set $s$ refines $\omega_0$, then $p\in F_0(x)$
for any $x\in s$ and therefore $f_0$ is a selection of $F_0|_{\Sigma_0}$.

Suppose that we already constructed $f_n$ --- a continuous selection of
$F_n|_{\Sigma_n}$. Let us define $f_{n+1}$ on arbitrary element $Z$
of discrete collection $\sigma_{n+1}$.
Since $\Sigma$ is locally finite, the set $A=Z\cap \Sigma_n$
is closed in $X$.
Since $f_n$ is a selection of $F_n$, then $f_n(A)$ is contained in $F_n(Z)$.
Since the pair $F_n(Z)\subset \cap_{x\in Z} F_{n+1}(x)$
is $X$-connected, we can extend $f_n|_A$ to a mapping
$f_n'\colon Z\to \cap_{x\in Z} F_{n+1}(x)$.
Clearly, $f_n'$ is a selection of $F_{n+1}|_Z$.
We define $f_{n+1}$ on the set $Z$ as $f_n'$.

Finally, we define $f$ to be equal to $f_n$ on the set $\Sigma_n$.
\end{proof}

\begin{thm} \label{thmusp}
Let $L$ be a finite $CW$-complex and
$F\colon X\to Y$ be a multivalued mapping of paracompact C-space
$X$ of extension dimension $\ed X\le [L]$ to a topological space $Y$.
If $F$ admits infinite fiberwise $[L]_c$-connected filtration of
strongly l.s.c. multivalued mappings,
then $F$ has a singlevalued continuous selection.
\end{thm}

\begin{proof}
By Theorem~\ref{thmtop}, the mapping $F'\colon X\to Y\times Q$
defined as $F'(x)=F(x)\times Q$ contains a stably $[L]$-connected
filtration of multivalued mappings.
By Theorem~\ref{thmsel} $F'$ has a singlevalued continuous selection $f'$.
Then the mapping $f=\pr_Y\circ f'$ is
a singlevalued continuous selection of $F$.
\end{proof}

\section{Hurewicz theorem}\label{S:Hurewicz}

The proof of the following theorem is similar to the proof
of Theorem~2.4 from~\cite{ChV}.

\begin{thm} \label{thmHur}
Let $f\colon X\to Y$ be a closed mapping of $k$-space $X$ onto
paracompact $C$-space $Y$.
Suppose that $\ed Y\le [M]$ for a finite $CW$-complex $M$.
If for every point $y\in Y$ and for every compactum $Z$
with $\ed Z\le [M]$ we have $\ed(f^{-1}(y)\times Z)\le [L]$
for some $CW$-complex $L$, then $\ed X\le [L]$.
\end{thm}

\begin{proof}
Suppose $A\subset X$ is closed and $g\colon A\to L$ is a map.
We are going to find a continuous extension $\wt g\colon X\to L$ of $g$.
Let $K$ be the cone over $L$ with a vertex $v$.
We denote by $C(X,K)$ the space of all continuous maps
from $X$ to $K$ equipped with the compact-open topology.
We define a multivalued map $F\colon Y\to C(X,K)$ as follows:
$$  F(y)=\{h\in C(X,K)\mid h(f^{-1}(y))\subset K\setminus \{v\}
    \text{ and } h|_A=g \}.  $$

{\it Claim.} $F$ admits continuous singlevalued selection.

\noindent If $\varphi\colon Y\to C(X,K)$ is a continuous selection
for $F$, then the mapping $h\colon X\to K$ defined
by $h(x)=\varphi(f(x))(x)$ is continuous on every compact subset
of $X$ and because $X$ is a $k$-space, $h$ is continuous.
Since $\varphi(f(x))\in F(f(x))$ for every $x\in X$,
we have $h(X)\subset K\setminus \{v\}$.
Now if $\pi\colon K\setminus \{v\}\to L$ denotes the
natural retraction, then $\wt g =\pi\circ h\colon X\to L$
is the desired continuous extension of $h$.

{\it Proof of the claim.}
We are going to apply Theorem~\ref{thmusp} to infinite filtration
$F\subset F\subset F\subset\dots$.
To do this, we have to show that $F$ is strongly l.s.c.
and that the pair $F(y)\subset F(y)$ is $[M]_c$-connected
for every point $y\in Y$.

First, we show that $F$ is strongly l.s.c.
Let $y_0\in Y$ and $P\subset F(y_0)$ be compact.
We have to find a neighborhood $V$ of $y_0$ in $Y$
such that $P\subset F(y)$ for every $y\in V$.
For every $x\in X$ define a subset $P(x)=\{h(x)\mid h\in P\}$ of $K$.
Since $P\subset C(X,K)$ is compact and $X$ is a $k$-space,
by the Ascoli theorem, each $P(x)$ is compact and $P$ is evenly continuous.
This easily implies that the set
$W=\{x\in X\mid P(x)\subset K \setminus \{v\}\}$
is open in $X$ and, obviously, $f^{-1}(y_0)\subset W$.
Since $f$ is closed, there exists a neighborhood $V$ of $y_0$ in $Y$
with $f^{-1}(V)\subset W$.
Then, according to the choice of $W$ and the definition of $F$,
we have $P\subset F(y)$ for every $y\in V$.

Fix an arbitrary point $y\in Y$. Let us prove that the pair
$F(y)\subset F(y)$ is $[M]_c$-connected.
Consider a pair of compacta $B\subset Z$ where $\ed Z\le [M]$
and a mapping $\varphi\colon B\to F(y)$.
Since $B\times X$ is a $k$-space (as a product of a compact
space and a $k$-space), the map $\psi\colon B\times X\to K$
defined as $\psi(b,x)=\varphi(b)(x)$ is continuous.
Extend $\psi$ to a set $Z\times A$ letting $\psi(z,a)=g(a)$.
Clearly, $\psi$ takes the set
$Z\times f^{-1}(y)\cap (Z\times A\cup B\times X)$
into $K\setminus \{v\}\cong L\times [0,1)$.
Since $\ed(Z\times f^{-1}(y))\le [L]$, we can extend $\psi$ over the set
$Z\times f^{-1}(y)$ to take it into $K\setminus \{v\}$.
Finally extend $\psi$ over $Z\times X$ as a mapping into $AE$-space $K$.
Now define an extension $\wt\varphi\colon Z\to F(y)$ of the mapping $\varphi$
by the formula $\wt\varphi(z)(x)=\psi(z,x)$.
\end{proof}

\begin{cor}[cf. Theorem~2.25 from \cite{DRS}]
Let $f\colon X\to Y$ be a mapping of finite-dimensional compacta
where $\ed Y=[M]$ for finite $CW$-complex $M$.
If for some $CW$-complex $L$ we have $\ed(f^{-1}(y)\times Y)\le [L]$
for every point $y\in Y$, then $\ed X\le [L]$.
\end{cor}

\begin{proof}
By Theorem~6.3 from \cite{DD} for any compactum $Z$ with
$\ed Z\le \ed Y$ we have $\ed(f^{-1}(y)\times Z)\le [L]$.
Thus, we can apply Theorem~\ref{thmHur}
\end{proof}

\section{Acknowledgments}

Authors wish to express their indebtedness to S.M.~Ageev for outlining
the proof of Theorem~\ref{thmtop} and to A.N.~Dranishnikov
for helpful discussions during the development of this work.

\end{document}